\documentclass[12pt]{article}

\pdfoutput=1
\usepackage{amsfonts}
\usepackage{amssymb}
\usepackage{amsmath}
\usepackage{amsthm}
\usepackage{wasysym}
\usepackage{graphicx}
\usepackage{geometry}
\usepackage{color}

\newcommand{\real}{\mathbb{R}}
\newcommand{\comp}{\mathbb{C}}
\newcommand{\nat}{\mathbb{N}}

\newcommand{\N}{\mathbb{N}}
\newcommand{\intd}{\mathrm{d}}
\newcommand{\epsi}{\varepsilon}
\newcommand{\vek}[1]{\textbf{#1}}

\newcommand{\norm}[2]{\left\|#1\right\|_{#2}}
\newcommand{\normLone}[1]{\left\|#1\right\|_{L^1}}
\newcommand{\normLtwo}[1]{\left\|#1\right\|_{L^2}}

\newcommand{\supp}{\mathop{\rm supp}}
\newcommand{\nnz}{\mathop{\rm nnz}}
\newcommand{\vj}{\,{\vek{j}}}
\newcommand{\bI}{{\vek{I}}}
\newcommand{\bk}{{\vek{k}}}

\newtheoremstyle{custom}{3pt}{3pt}{}{}{\bfseries}{:}{.5em}{}
\theoremstyle{custom}
\newtheorem{beispiel}  {Example}[section]
\newtheorem{definition}[beispiel]{Definition}
\newtheorem{satz}      [beispiel]{Theorem}    
\newtheorem{lemma}     [beispiel]{Lemma}

\newtheorem{bem} 		[beispiel]{Remark}
\newtheorem{proposition}[beispiel]{Proposition}

\numberwithin{equation}{section}

\begin{document}


\title{Discretization of transfer operators using a sparse hierarchical tensor basis -- the Sparse Ulam method}
\author{Oliver~Junge and P\'eter~Koltai\\Faculty for Mathematics\\Technische Universit\"at M\"unchen}

\date{\today}
\maketitle
\begin{abstract}
The global macroscopic behaviour of a dynamical system is encoded in the eigenfunctions of a certain transfer operator  associated to it.  For systems with low dimensional long term dynamics, efficient techniques exist for a numerical approximation of the most important eigenfunctions, cf.\ \cite{DeJu99a}.  They are based on a projection of the operator onto a space of piecewise constant functions supported on a neighborhood of the attractor -- Ulam's method.  

In this paper we develop a numerical technique which makes Ulam's approach applicable to systems with higher dimensional long term dynamics.  It is based on ideas for the treatment of higher dimensional partial differential equations using \emph{sparse grids} \cite{Ze91a,BuGr04a}.  We develop the technique, establish statements about its complexity and convergence and present two numerical examples.
\end{abstract}

\section{Introduction}

Recently, numerical techniques have been developed which enable a coarse grained, yet global statistical analysis of the long term behaviour of certain dynamical systems.  The basic algorithmic approach is to construct a box covering of some set of interest in phase space (e.g.\ the attractor of the system) \cite{DeHo96a,DeHo97a}.  The cells in this covering then constitute the states of a finite Markov chain.  The transition matrix of this chain (i.e.\ the matrix of transition probabilities between the boxes) can be viewed as a finite approximation to the transfer (or Frobenius-Perron) operator of the system.  This operator describes how probability distributions on phase space evolve according to the dynamical system under consideration.  In certain cases and in the appropriate functional analytic setting, eigenmodes of this operater can be used to charaterize the long term behaviour of the dynamics.  Certain stationary distributions of the operator characterize how frequently typical trajectories visit certain parts of phase space.  Eigenmodes at roots of unity enable the detection of macroscopic cycles in the dynamics and eigenmodes at real eigenvalues close to one yield a decomposition of phase space into \emph{almost invariant sets}, i.e.\ sets for which the probability for a typical point to be mapped back into the set is large \cite{DeJu99a}.  The latter concept has e.g.\ been used in order to detect and compute biomolecular conformations, cf. \cite{DeDeJu99a,Sc99a,ScHuDe01a,ScHu03a,DeSc04a}.

Formally, the construction of the Markov chain can be viewed as projecting the transfer operator onto the space of functions which are piecewise constant on the elements of the box covering.  Ulam conjectured \cite{Ul60a} that for maps on the interval, the stationary distribution of the chain converges to an invariant density (i.e.\ a stationary distribution) of the map.  This has been proved for certain expanding maps by Li \cite{Li76a} and since then for various special classes of maps or stochastic processes also in higher dimensions \cite{DiLi91a,DiDuLi93a,Fr95a,DiZh96a,Fr98a,DeJu99a}.

Ulam's method in combination with the subdivision approach from \cite{DeHo96a,DeHo97a} for the computation of the box covering  works fine for systems with a low dimensional attractor, cf. also \cite{DeFrJu01a,DeJu02b}. For systems with higher dimensional long term dynamics the approach becomes inefficient due to the \emph{curse of dimension}: the number of boxes in the covering scales exponentially in the dimension of the attractor.  Adaptive approaches to the construction of the box covering \cite{DeJu98a,Ju01a} do not remedy this fact.  

In this paper we propose to attack this discretization task using ideas from \emph{sparse grids} \cite{Sm63a,Ze91a,BuGr04a}.  In this approach, which is e.g.\ being used in the numerical solution of partial differential equations on higher dimensional domains, a basis of $[0,1]^d$ is build from a  hierarchical basis of $[0,1]$ via a tensor product construction.  The entire basis can be decomposed into subspaces which are spanned by basis functions of the same level of the 1d hierarchy in each factor.  To each subspace one can associate its approximation \emph{benefit} and its \emph{cost} (which is typically given by its dimension).  The idea of the sparse grid approach is to assemble a finite dimensional approximation space by choosing only those subspaces whith the highest benefit to cost ratio.  

In order to discretize the Frobenius-Perron operator, we employ a piecewise constant sparse hierarchical tensor basis (i.e.\ using the Haar system as the underlying 1d basis).  This basis provides an approximation error of $\mathcal{O}(n^{-1}\cdot(\log\sqrt{n})^{d-1})$ for functions with bounded first derivatives, requiring a computational effort of $\mathcal{O}(n\cdot(\log\sqrt{n})^{d-1})$ (where $n$ denotes the number of degrees of freedom in one coordinate direction and $d$ is the dimension of phase space).  In comparison, the standard Ulam basis requires $\mathcal{O}(n^d)$ basis functions in order to obtain an approximation error of $\mathcal{O}(n^{-1})$.

The paper is structured as follows: in Section~\ref{sec:FP} we collect relevant basic concepts from dynamical systems theory, in particular Ulam's method.  In Section~\ref{sec:SP} we develop the Sparse Ulam method by constructing the hierarchical tensor basis, deriving approximation properties, outlining the construction of the optimal approximation subspace and comparing cost and accuracy of the new method with the standard Ulam approach.  The section closes with statements about the convergence properties.  In Section~\ref{sec:comp} we collect considerations concerning an efficient implementation of our approach.  In particular, we derive estimates on the computational effort as a function of the required accuracy.  Section~\ref{sec:examples} presents two numerical examples: a comparison with Ulam's method for a three dimensional map with a smooth invariant density and a computation of the leading eigenfunctions of the transfer operator for a four-dimensional map, constructed via a tensor product from two two-dimensional standard maps.

Our implementation of the Sparse Ulam method as well as the code for the example computations is freely available from the homepage of the authors.

\section{Transfer operators and Ulam's method}
\label{sec:FPO_ulam}

\subsection{Long term dynamics and the Frobenius-Perron operator}\label{sec:FP}

 Let $S: X\rightarrow X$, $ X\subset\real^d$, be a discrete dynamical system which is measureable w.r.t the Borel-$\sigma$-algebra $\mathcal{B}$ on $ X$. Let $\mathcal{M}_\comp$ be the set of all bounded complex valued measures on $(X,\mathcal{B})$ and $\mathcal{M}\subset\mathcal{M}_\comp$ be the subset of probability measures. The \emph{Frobenius-Perron operator} (or \emph{transfer operator}) $P:\mathcal{M}_\comp\rightarrow\mathcal{M}_\comp$,
\begin{equation}
	P\mu = \mu\circ S^{-1},
	\label{eq:FPO_gen}
\end{equation}
describes how (probability) measures on phase space evolve according to the dynamics defined by $S$.
A measure is called \emph{invariant} if it is a fixed point of $P$. A set $A\subset X$ is called invariant if $A=S^{-1}(A)$.  An invariant probability measure $\mu$ is \emph{ergodic} if every invariant set has either full or zero $\mu$-measure.  Birkhoff's \emph{ergodic theorem} \cite{Bi31a} states that ergodic measures characterize the long time behaviour of the system: Let $\mu$ be ergodic and $\varphi: X\to\real$ be a $\mu$-integrable \emph{observable}, then
\begin{equation}\label{eq:ergodic theorem}
\lim_{n\rightarrow\infty}\frac{1}{n}\sum^{n-1}_{k=0}\varphi(S^{k}(x))=\int_{ X}\varphi\:\intd \mu
\end{equation}
for $\mu$-almost all $x\in X$.  

\begin{definition}	\label{def:SRB}
	A probability measure $\mu$ is called \emph{SRB measure} or \emph{natural invariant measure} if (\ref{eq:ergodic theorem}) holds for continuous observables $\varphi$ and all points $x\in U$ in a set $U\subset X$ with positive Lebesgue measure.	\end{definition}
SRB measures are defined via a property which we would like them to have. But how does one see whether a measure is SRB? After all, equation \eqref{eq:ergodic theorem} is not easy to check in general. On the other hand, if $\mu$ is an ergodic measure which is absolutely continuous w.r.t.\ the Lebesgue-measure $m$, i.e.\ if there is a \textit{density} $f$ with
$ \mu(A) = \int_{A}f\ \intd m$ for all measurable $A$,
then $\mu$ is SRB.

Using \eqref{eq:FPO_gen}, we can directly define the Frobenius-Perron operator on Lebesgue integrable functions $f: X\to\comp$:
\begin{equation}
	\int_A Pf\ \intd m = \int_{S^{-1}(A)} f\ \intd m \qquad \forall A\in\mathcal{B}.
	\label{eq:FPO_density}
\end{equation}
If $S$ is differentiable, we obtain the explicit expression
\[ Pf(x) = \sum_{y\in S^{-1}(x)}\frac{f(y)}{\left|DS(y)\right|}. \]

\subsection{Almost invariance}\label{ssec:almost}

Invariant measures (or densities) are fixed points of the Frobenius-Perron operator, i.e.\ eigenmeasures resp.\ -functions at the eigenvalue 1.  Eigenvectors at eigenvalues close to one are related to \emph{almost invariant sets}:
Intuitively, an \emph{almost invariant set} of $S$ is a subset $A\subset  X$ such that the \emph{invariance ratio}
\[
\rho_m(A) := \frac{m(S^{-1}(A)\cap A)}{m(A)}
\]
is close to $1$, i.e.\ a point which is chosen randomly from $A$ with respect to the measure $m$ maps into $A$ with high probability.   More precisely, we say that
\begin{definition}
	A subset $A\subset X$ is \emph{$\varrho$-almost-invariant} w.r.t.\ the probability measure $\mu$ if $\mu(A)\neq 0$ and
	\begin{equation}\label{eq:almost_inv}
		\rho_\mu(A) = \varrho.
	\end{equation}
\end{definition}
Let $\nu\in\mathcal{M}_\comp$ be an eigenmeasure of $P$ at an eigenvalue $\lambda \neq 1$. Since
$\lambda\nu( X) = (P\nu)( X) = \nu\left(S^{-1}\left( X\right)\right) = \nu( X)$, it follows that $\nu( X)=0$. In particular, if $\lambda< 1$ and $\nu$ are real, then there are two positive real measures $\nu^+, \nu^-$ with disjoint supports  such that $\nu=\nu^+ - \nu^-$ (\emph{Hahn-Jordan decomposition}).  The following theorem relates the invariance ratios of the supports of $\nu^+$ and $\nu^-$ to the eigenvalue $\lambda$.
\begin{satz}\cite{DeJu99a}
	Let $\nu$ be a normalized\footnote{i.e. $|\nu|( X)=1$.} eigenmeasure of $P$ at the real eigenvalue $\lambda < 1$. Then 	\begin{equation}
		\rho_{|\nu|}(A^+)+\rho_{|\nu|}(A^-) = \lambda+1,
		\label{eq:eigenval_inv}
	\end{equation}
	where $A^+=\supp \nu^+$ and $A^-=\supp \nu^-$.
\end{satz}

\subsection{Ulam's method}\label{subsec:Ulam}

In order to approximate the (most important) eigenfunctions of the Frobenius-Perron operator, we have to discretize the corresponding infinite dimensional eigenproblem.
Ulam \cite{Ul60a} proposed to project the $L^1$ eigenvalue problem $Pf=\lambda f$ into a finite dimensional subspace of piecewise constant functions: 
Let $(V_n)_{n\in\nat}$ be a sequence of \textit{approximation subspaces} of $L^1$ with $\dim V_n = n$ and let $Q_n:L^1\rightarrow V_n$ be corresponding projections into $V_n$. The sequences $(V_n)$ and $(Q_n)$ should be chosen such that $Q_{n}$ converges pointwise to the identity on $L^1$.  We define the \emph{discretised Frobenius-Perron operator} as 
\[
P_n := Q_nP\mid_{V_n}.
\]  

We now choose the approximation spaces to be spanned by piecewise constant functions. To this end, let $\mathcal X_n = \{X_1,\ldots,X_n\}$ be a disjoint partition of $X$ with $m(X_i)\to 0$ as $n\to\infty$.
 Define $V_n := \text{span}\{\chi_{1},\ldots,\chi_{n}\}$,
where $\chi_{i}$ denotes the characteristic function of $X_i$.\\
Further, let
\[ Q_{n}h:=\sum_{i=1}^{n}c_i\chi_{i}\qquad\text{ with } \qquad c_i:=\frac{1}{m(I_i)}\int_{I_i}h\:\intd m, \]
yielding $P_{n}\Delta_{n}^{+}\subseteq\Delta_{n}^{+}$ and $P_{n}\Delta_{n}\subseteq\Delta_{n}$, where
$	\Delta_{n}  :=  \left\{h\in V_n : \int |h|\; dm = 1\right\}$ and $\Delta_{n}^{+}  :=  \left\{ h\in\Delta_n : h\geq 0\right\}.
$
Due to Brouwer's fixed point theorem there always exists an approximative invariant density $f_n = P_nf_n\in\Delta_n^+$.  The matrix representation of the linear map $P_n:\Delta_n\rightarrow\Delta_n$ w.r.t.\ the basis of characteristic functions is given by the \textit{transition matrix} with entries
\begin{equation}
	p_{ij}=\frac{m(X_j\cap S^{-1}(X_i))}{m(X_j)}.
	\label{eq:trans_rate}
\end{equation}
Ulam conjectured \cite{Ul60a} that if $P$ has a unique stationary density $f$, then a sequence $(f_n)_{n\in\nat}$ converges to $f$ in $L^1$.  It is still an open question under which conditions on $S$ this is true in general.  Li \cite{Li76a} proved the conjecture for expanding, piecewise continuous interval maps, Ding an Zhou  \cite{DiZh96a} for the corresponding multidimensional case.  

In \cite{DeJu99a}, Ulam's method was applied to a small random perturbation of $S$ which might be chosen such that the corresponding transfer operator is compact on $L^2$.  In this case, perturbation results \cite{Ka84a} (section IV.3.5.) for the spectrum of compact operators imply convergence.


\subsection{Computing the transition matrix}	\label{ssec:comp_trm}

The computation of one matrix entry \eqref{eq:trans_rate} requires a $d$-dimensional quadrature.  A standard approach to this is Monte-Carlo quadrature (also cf. \cite{Hu94a}), i.e.
\begin{equation}
	p_{ij} \approx \frac{1}{K}\sum_{k=1}^{K}\chi_i\left(S(x_k)\right),
	\label{eq:computePij}
\end{equation}
where the points $x_1,\ldots,x_K$ are chosen i.i.d from $X_j$ according to a uniform distribution. In   \cite{GuDeKr97a}, a recursive exhaustion technique has been developed in order to compute the entries to a prescribed accuracy. However, this approach relies on the availability of local Lipschitz estimates on $S$ which might not be cheaply computable in the case that $S$ is given as the time-$T$-map of a differential equation.

For the Monte-Carlo technique, consider a uniform partition of the unit cube into $M^{d}$ congruent cubes of edge length $1/M$. Let $P_M$ denote the transition matrix for this partition and let $\tilde{P_M}$ be its Monte-Carlo approximation.  
According to the central limit theorem (and its error-estimate, the Berry-Ess\'een theorem \cite{Fe71a}) we have\footnote{We write $a(K)\apprle b(K)$ if there is a  constant $c>0$ independet of $K$ such that $a(K) \leq cb(K)$.}
\begin{equation}
	|\tilde{p}_{ij}-p_{ij}| \apprle 1/{\sqrt{K}}
	\label{eq:entr_err_ulam}
\end{equation}
for the absolute error of the entries of $\tilde P_M$.
As a consequence, we need 
\begin{equation}
	\varrho \apprge \frac{M^{d}}{TOL^2}
	\label{eq:point_nr_ulam}
\end{equation}
sample points in total in order to achieve an absolute error of less than $TOL$ for all the entries $p_{ij}$.
Note that the accuracy of the entries of $\tilde P_M$ imposes a restriction on the achievable accuracy of the eigenvectors of $P_M$. 

\section{The Sparse Ulam method}\label{sec:SP}

A naive application of Ulam's method to higher dimensional systems suffers from the \emph{curse of dimension}: in order to achieve an $L^1$-accuracy of $\mathcal O(\varepsilon)$ one needs an approximation space of dimension $\mathcal{O}\left(\epsi^{-d}\right)$ -- translating into a prohibitively large computational effort for higher dimensional systems.  There is a remedy to this problem for systems with low dimensional long term dynamics \cite{DeHo97a,DeJu99a}: the idea is to first compute a covering of the attractor of the system.  On this (low dimensional) covering, Ulam's method can successfully be applied. 

To avoid the exponential growth of complexity in the system (or attractor) dimension, we now follow an idea which was originally developed for quadrature problems \cite{Sm63a} and used for the treatment of higher dimensional partial differential equations, cf. for example \cite{Ze91a,BuGr04a}: \emph{sparse grids}.  In fact, we change from the standard Ulam basis to a sparse hierarchical one in order to obtain a better cost/accuracy relation. In the following, we discuss the chosen basis in detail, as well its advantages and disadvantages.

\subsection{The Haar basis}

We describe the Haar basis on the $d$-dimensional unit cube $[0,1]^{d}$, deriving the multidimensional basis functions from the one dimensional ones, see e.g.\ \cite{GrOsSc99a}. Let
\begin{equation}
	f_\text{Haar}(x)= -\text{sign}(x)\cdot(\left|x\right|\leq 1), 		\label{eq:base}
\end{equation}
where $(\left|x\right|\leq 1)$ equals $1$, if the inequality is true, otherwise 0.  A basis function of the Haar basis is defined by the two parameters \emph{level} $i$ and  \emph{center (point)} $j$:
\begin{equation}
	f_{i,j}(x):= \left\{\begin{array}{ll}
						1 & \text{if } i=0, \\
						2^{\frac{i-1}{2}}\cdot f_\text{Haar}\left(2^{i}\left(x-x_{i,j}\right)\right) & \text{if } i\geq 1,
				 \end{array}\right.		\label{eq:1Dbasisfun}
\end{equation}
where
\begin{equation}
	x_{i,j}:=(2j+1)/2^{i}, \qquad j\in\{0,\ldots,2^{i-1}-1\}.
\end{equation}
A $d$-dimensional basis function is constructed from the one dimensional ones using a tensor product construction:
\begin{equation}
	\varphi_{\ell,\vj}(x) := \prod_{i=1}^{d}f_{\ell_{i},\vj_{i}}(x_i),
	\label{eq:dDbasisfun}
\end{equation}
for $x=(x_1,\ldots,x_d)\in [0,1]^d$. Here $\ell=(\ell_1,\ldots,\ell_d)$, $\ell_i\in\{0,1,2,\ldots\}$, denotes the level of the basis function and $\vj=(\vj_1,\ldots,\vj_d)$, $\vj_{i}\in \{ 0,\ldots,2^{\ell_{i}}-1\}$, its center.
\begin{satz}[Haar basis]\label{thrm:haar}
	The set
	\[ H=\left\{f_{i,j}\mid i\in\nat_{0}, j\in\{0,\ldots,2^{i}-1\} \right\} \]
	is an orthonormal basis of $L^2([0,1])$, the \emph{Haar basis}.	Similarly, the set
	\begin{equation*}
	H^d = \left\{\varphi_{\ell,\vj}\mid \ell\in\nat_{0}^d, \vj_i\in\{0,\ldots,2^{\ell_i}-1\} \right\}
	\end{equation*}
	is an orthonormal basis of $L^2([0,1]^d)$.
\end{satz}

\begin{figure}[htb]
\centering
\includegraphics[width=1.00\textwidth]{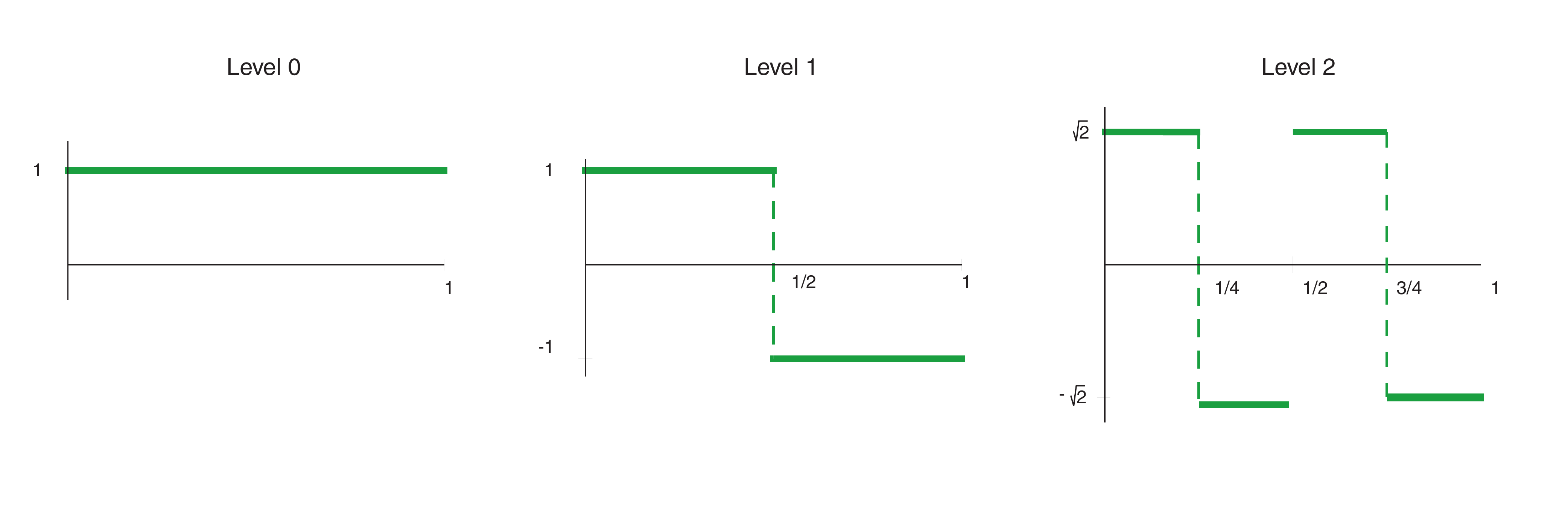}
\caption{First three levels of the 1D Haar basis}
\label{fig:haar1D}
\end{figure}

Figure \ref{fig:haar1D} shows the basis functions of the first three levels of the one dimensional Haar basis. 
It will prove useful to collect all basis functions of one level in one subspace:
\begin{equation}
	W_{\ell} := \text{span}\left\{ \varphi_{\ell,\vj}\mid \vj_{i}\in \{0,\ldots,2^{\ell_i}-1\} \right\}, \qquad \ell\in\nat_{0}^{d}.
	\label{eq:Hierarchie}
\end{equation}
Consequently, $L^2=L^2([0,1]^d)$ can be written as the infinite direct sum of the subspaces $W_{\ell}$,
\begin{equation}
	L^{2} = \bigoplus_{\ell\in\nat_{0}^{d}}W_{\ell}.
	\label{eq:directsum}
\end{equation}
In fact, it can also be shown that $L^1 = L^1([0,1]^d)=\bigoplus_{\ell\in\nat_{0}^{d}}W_{\ell}$ as well. To see this, note that $\bigoplus_{\ell\in\vek{I}^{d}}W_{\ell}$ with $\vek{I} = \{\ell \mid \norm{\ell}{\infty}\leq n\}$ is the space of characteristic functions supported on the uniform decomposition of the unit cube in $2^n$ subcubes in every direction.
Moreover, we have
\begin{equation}
	\text{dim }W_{\ell} = \prod_{i=1}^{d}2^{\max\{0, \ell_{i}-1\}} = 2^{\sum_{\ell_i\neq0}\ell_{i}-1}.
	\label{eq:dimWl}
\end{equation}
In order to get a finite dimensional approximation space most appropriate for our purposes, we are going the choose an optimal finite subset of the basis functions $\varphi_{\ell,\vj}$. Since in general we do not have any a priori information about the function to be approximated,   and since all basis functions in one subspace $W_{\ell}$ deliver the same contribution to the approximation error we will use either all or none of them. In other words, the choice for the approximation space is transferred to the level of subspaces $W_{\ell}$. 

\subsection{Approximation properties}

The choice of the optimal set of subspaces $W_\ell$ relies in the contribution of each of these to the approximation error. The following statements give estimates on this.

\begin{lemma}
	Let $f\in C^{1}([0,1])$ and let $c_{i,j}$ be its coefficients with respect to the Haar basis, i.e.\ $f=\sum_{ij}c_{i,j}f_{i,j}$. Then for $i>0$ and all $j$
	\[ 
	\left|c_{i,j}\right|\leq 2^{-\frac{3i+1}{2}}\|f'\|_\infty.
	\]
	For $f\in C^{1}\left([0,1]^d\right)$ we analogously have for $\ell\neq 0$ and all $\vj$
	\[ 
	\left|c_{\ell,\vj}\right|\leq 2^{-\left(\sum_{\ell_{i}\neq 0}3\ell_{i}+1\right)/2}\prod_{\ell_{i}\neq0} \|\partial_i f\|_{\infty}. 
	\]
\end{lemma}
\begin{proof}
For $i\geq 1$
\begin{eqnarray*}
2^{\frac{1-i}{2}}c_{ij} & = & \int_{x_j-2^{-i}}^{x_j}f-\int^{x_j+2^{-i}}_{x_j}f \\
& = & \int_{x_j-2^{-i}}^{x_j} \left(f(x_j)+\int_{x_j}^{x}f'\right)\intd x - \int_{x_j}^{x_j+2^{-i}} \left(f(x_j)+\int_{x_j}^{x}f'\right)\intd x
	\end{eqnarray*}
	and thus
\[ 
2^{\frac{1-i}{2}}\left|c_{ij}\right| \leq 2\|f'\|_\infty \int_{0}^{2^{-i}}x\, \intd x,
\]
which yields the claimed estimate for the 1D case.  The bound in the $d$-dimensional case follows similarly.
\end{proof}

Using this bound on the contribution of a single basis function to the approximation of a given function $f$, we can derive a bound on the total contribution of a subspace $W_\ell$. For $f_\ell\in W_\ell$ 
\begin{eqnarray}
	\normLone{f_{\ell}} & \leq & 2^{-\sum_{\ell_i\neq0}(\ell_i+1)}\prod_{\ell_{i}\neq0} \|\partial_i f\|_\infty,
	\label{eq:L1beitrag} \\
	\normLtwo{f_{\ell}} & \leq & 2^{-\sum_{\ell_i\neq 0}(\ell_i+3)/2}\prod_{\ell_{i}\neq 0} \|\partial_i f\|_\infty.
	\label{eq:L2beitrag}
\end{eqnarray}

\subsection{The optimal subspace}

The main idea of the sparse grid approach is to choose \emph{cost} and (approximation) \emph{benefit} of the approximation subspace in an optimal way.  We briefly sketch this idea here, for a detailed exposition see \cite{Ze91a,BuGr04a}.  For a set ${\bf I}\subset \N^d_0$ of multiindices we define
\[
W_{\bf I} = \bigoplus_{\ell\in\bf I} W_\ell.
\]
Correspondingly, for $f\in L^1$, let $f_{\bf I}=\sum_{\ell\in\bf I} f_\ell$, where $f_\ell$ is the orthogonal projection of $f$ onto $W_\ell$.  We define the \emph{cost} $C(\ell)$ of a subspace $W_\ell$ as its dimension,
\[
C(\ell)=\dim W_\ell = 2^{\sum_{\ell_i\neq 0}\ell_i-1}.
\]

Since
\begin{eqnarray}
	\left\|f-f_{\vek{I}}\right\| & \leq & \sum_{\ell\notin\vek{I}}\left\|f_{\ell}\right\| \label{eq:accuracy1} 
															  =  \sum_{\ell\in\nat_{0}^{d}}\left\|f_{\ell}\right\| - \sum_{\ell\in\vek{I}}\left\|f_{\ell}\right\|,	\label{eq:accuracy2}
\end{eqnarray}
the guaranteed increase in accuracy is bounded by the contribution of a subspace $W_{\ell}$ which we add to the approximation space.  We therefore define the \emph{benefit} $B(\ell)$ of $W_\ell$ as the upper bound on its $L_1$-contribution as derived above,
\begin{equation}
	B(\ell) = 2^{-\sum_{\ell_i\neq0}(\ell_i+1)}.
	\label{eq:benefit}
\end{equation}
Note that we omited the factor involving derivatives of $f$. The reason is that it does not affect the solution of the optimization problem \eqref{eq:opt_problem}

\begin{figure}[htb]
  \centering
	\includegraphics[width=0.50\textwidth]{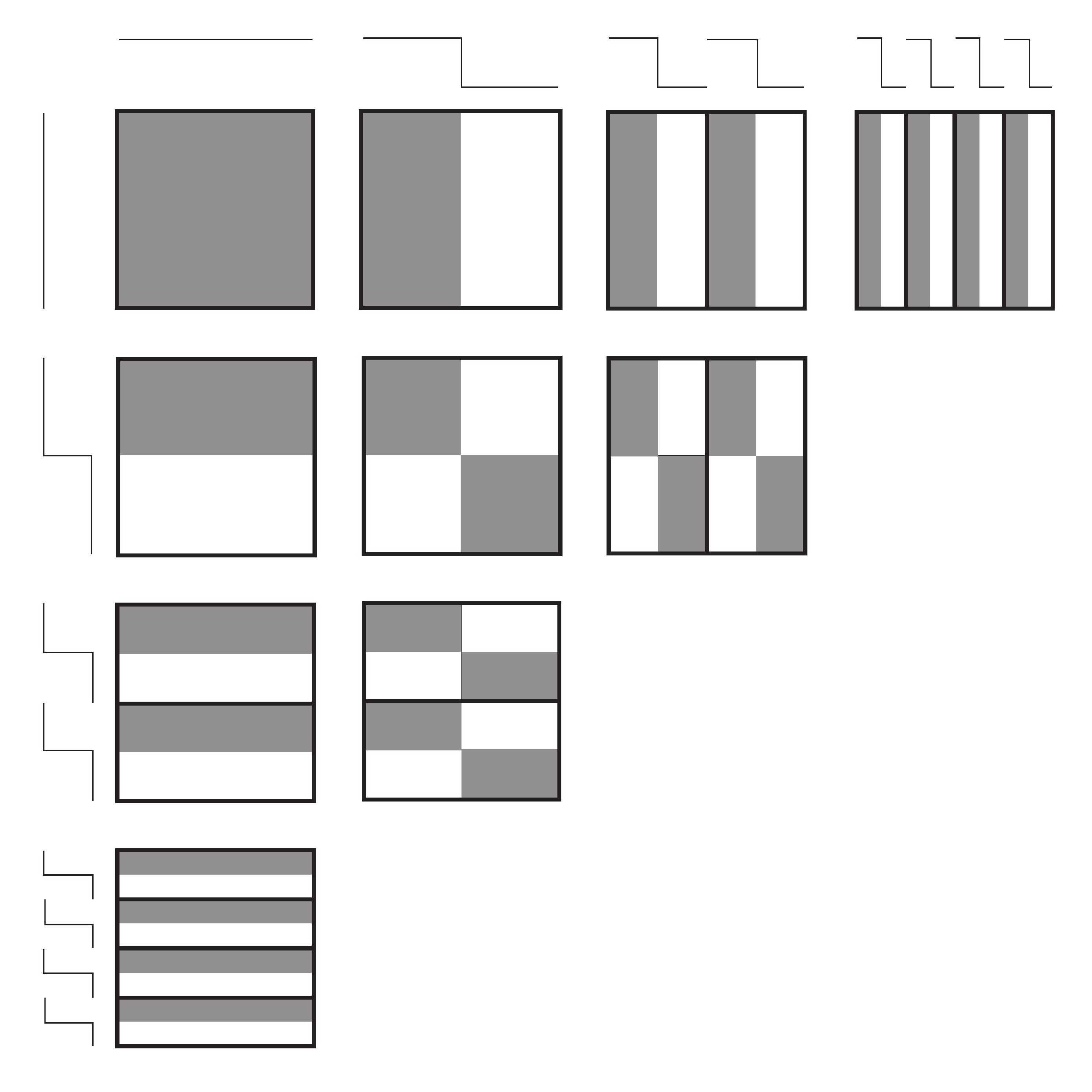}
	\caption{$3^{rd}$ level sparse basis in 2D.	Shaded means value 1, white means value $-1$, thicker lines are support boundaries.}
	\label{fig:Haar_basis2D_b_transparent}
\end{figure}

Let $C(\vek{I})=\sum_{\ell\in\vek{I}}C(\ell)$ and $B(\vek{I})=\sum_{\ell\in\vek{I}}B(\ell)$ be the total cost and the total benefit of the approximation space $W_{\bf I}$. In order to find the optimal approximation space we are now solving the following optimization problem:  Given a bound  $c>0$ on the total cost, find an approximation space $W_{\bf I}$ which solves
\begin{equation}
	\max_{C({\bI})\leq c}B({\bI}).
	\label{eq:opt_problem}
\end{equation}
One can show (cf.\ \cite{BuGr04a}) that ${\bf I}\subset\N^d_0$ is an optimal solution to (\ref{eq:opt_problem}) iff 
\begin{equation}
	\frac{C(\ell)}{B(\ell)} = const\qquad\text{for }\ell\in \partial\bI,
	\label{eq:opt_cond1}
\end{equation}
where the \emph{boundary} $\partial\bI$ is given by
$\partial\bI = \{ \ell\in \bI \mid \ell'\in \bI, \ell'\geq \ell\ \Rightarrow \ell'=\ell\}$\footnote{$\ell'\geq \ell$ is meant componentwise}.
Using the definitions for cost and benefit as introduced above, we obtain
\begin{equation}
	  \frac{C(\ell)}{B(\ell)}= \frac{2^{\sum_{\ell_i\neq0}(\ell_{i}-1)}}{2^{-\sum_{\ell_{i}\neq0}(\ell_i+1)}} = 2^{2\sum_{\ell_i\neq0}\ell_i} = 2^{2|\ell|},
	\label{eq:opt_cond2}
\end{equation}
where $|\ell|$ means the 1-norm of the vector $\ell$.\\
The optimality condition (\ref{eq:opt_cond1}) thus translates into the simple condition
\begin{equation}
	|\ell| = const \quad\text{for }\ell\in\partial\bI.
	\label{eq:L1sparse_select}
\end{equation}
As a result, the optimal approximation space is $W_{\bI(N)}$ with
\begin{equation}
 	\bI(N) = \left\{\ell\in\N^d_0 \mid |\ell|\leq N \right\},
	\label{eq:indexset}
\end{equation}
where the \emph{level} $N=N(c)\in\N$ is depending on the chosen cost bound $c$.  Figure~\ref{fig:Haar_basis2D_b_transparent} schematically shows the basis functions of the optimal subspace in $2D$ for $N=3$.
\begin{bem}
	Because of the orthogonality of the Haar-basis in $L^2$ one can take the squared contribution as the benefit in the $L^2$-case (resulting in equality in \eqref{eq:accuracy1}).  In this case we obtain the optimality condition
	\begin{equation}
		\sum_{\ell_i\neq0}(\ell_{i}+1) = const\quad\text{for } \ell\in\partial\bI
		\label{eq:L2sparse_select}
	\end{equation}
	and correspondingly $W_\bI$ with
	\begin{equation}
	 	\bI(N) = \left\{\ell\in\N^d_0 : \sum_{\ell_i\neq0}(\ell_{i}+1)\leq N\right\},
		\label{eq:indexset2}
	\end{equation}
	$N=N(c)$, as the optimal approximation space.
\end{bem}

\subsection{The discretized operator}

Having chosen the optimal approximation space $V_N=W_{\bI(N)}$ we now build the corresponding discretized Frobenius-Perron operator $P_N$. Since the sparse basis
\begin{equation}
	B_N := \left\{ \varphi_{\ell,\vj}\mid |\ell|\leq N, \vj_{i}\in \{0,\ldots,2^{\ell_i}-1\}  \right\}
	\label{eq:Hierarchie1}
\end{equation}
is an $L^2$-orthogonal basis of $V_N$, the natural projection $Q_N:L^2\to V_N$ is given by
\begin{equation}
	Q_N f = \sum_{\varphi\in B_N}\left(\int f\varphi\right)\varphi.
	\label{eq:L2proj}
\end{equation}
All basis functions $\varphi\in B_N$ are piecewise constant and have compact support, so $Q_N$ is well defined on $L^1$ as well. Choosing an arbitrary enumeration, the (\emph{transition}) matrix of the discretized Frobenius-Perron operator
\[
P_N = Q_N\circ P
\]
with respect to $B_N$ has entries
\begin{equation}
	p_{ij} = \int \varphi_i\ P\varphi_j.
	\label{eq:FPO_sparse1}
\end{equation}
Writing $\varphi_i = \varphi_i^+-\varphi_i^- = |\varphi_i|\cdot(\chi^{+}_i-\chi^{-}_i)$, where $|\varphi_i|$ is the (constant) absolute value of the function over its support and $\chi^{+}_i$ and $\chi^{-}_i$ are the characteristic functions on the supports of the positive and negative parts of $\varphi_i$, we obtain
\begin{equation}
	p_{ij} = |\varphi_i||\varphi_j|\left(\int \chi^{+}_i\ P\chi^{+}_j - \int \chi^{-}_i\ P\chi^{+}_j - \int \chi^{+}_i\ P\chi^{-}_j + \int \chi^{-}_i\ P\chi^{-}_j\right),
	\label{eq:FPO_sparse2}
\end{equation}
which is, by \eqref{eq:trans_rate}
\begin{equation}
	p_{ij} = |\varphi_i||\varphi_j| \sum \pm m\left(X^{\pm}_j\cap S^{-1}\left(X^{\pm}_i\right)\right),
	\label{eq:FPO_sparse3}
\end{equation}
where $X_i^\pm=\supp\varphi_i^\pm$ and we add the 4 summands like in \eqref{eq:FPO_sparse2}. These can be computed in the same way as presented in section \ref{sec:FPO_ulam}.
 
\begin{bem}
	We note that
	\begin{itemize}
		\item[(a)] if the $i^{th}$ basis function is the one corresponding to $\ell = (0,\ldots,0)$, then 		
		\[ 
		p_{ij} = \delta_{ij}. 
		\]
		\item[(b)] The entries of $P_N$ are bounded via
			\[ \left|p_{ij}\right| \leq \sqrt{\frac{m(X_j)}{m(X_i)}} \leq 2^{N/2}.\]
		\item[(c)] If $P_N x = \lambda x$ with $\lambda \neq 1$, then $x_i = 0$ if the $i^{th}$ basis function is the one corresponding to $\ell = (0,\ldots,0)$.
		This follows from
		\begin{equation}\label{eq:ev_orth}
		 x_i \stackrel{(a)}{=} (e_i^{\top}P_N)x = e_i \lambda x = \lambda x_i. 
		 \end{equation}
		It is straightforward to show that this property is shared by every Ulam type projection method with a constant function as element of the basis of the approximation space. This observation is useful for the reliable computation of an eigenvector at an eigenvalue close to one (since it is badly conditioned):  (\ref{eq:ev_orth}) allows us to reduce the eigenproblem to the subspace orthogonal to the constant function.
	\end{itemize}
	With the given change in (c) are properties (a)-(c) valid for the numerical realisation as well.
\end{bem}

\subsection{Convergence}

As has been pointed out in the Introduction and in Section~\ref{subsec:Ulam}, statements about the convergence of Ulam's method exist in certain cases.  Note that for $N=kd$, $k=0,1,2,\ldots$, the approximation space $W_{\bI(N)}$ includes the Ulam approximation space $W_\ell$ with $\ell=(k,\ldots,k)$ and thus we obtain convergence of the Sparse Ulam method as a corollary to the convergence of Ulam's method in these cases from the following Lemma (which can be proved by standard arguments).
An open question is, if in general, the convergence of Ulam's method implies convergence of Sparse Ulam.

\begin{lemma}
	$\normLone{Q_N f-f}\stackrel{n\rightarrow\infty}{\longrightarrow}0$ for $f\in L^1$. 
\end{lemma}

\section{Complexity}\label{sec:comp}

In this section, we collect basic statements about the complexity of both methods.

\subsection{Cost and accuracy} \label{subsec:cost_accu}

We defined the total cost of an approximation space as its dimension and the accuracy via its \textit{contribution} or \textit{benefit}, see \eqref{eq:benefit}. In this section we derive a recurrence formula for these numbers, depending on the \textit{level} of the optimal subspaces and the system dimension. 

Let $C(N,d)$ be the dimension of $W_{\bI(N)}$ in phase space dimension $d$. Then
\begin{equation}
	C(N,d) = C(N,d-1) + \sum_{k=1}^{N}C(N-k,d-1)2^{k-1},
	\label{eq:rekurs_ges}
\end{equation}
since if $\ell = (*,\ldots,*,0)$, then the last dimension plays no role in the number of basis functions, and the total number of basis function's for such $\ell$'s is $C(N,d-1)$. If on the other hand $\ell = (*,\ldots,*,\ell_{d})$ with $\ell_d>0$, then the number of basis functions with such $\ell$'s is $C(N-\ell_d,d-1)2^{\ell_d-1}$, because there are $2^{\ell_d-1}$ one-dimensional basis functions of level $\ell_d$ possible for the tensor product in the last dimension.  For $d=1$ we simply deal with the standard Haar basis, so $C(N,1)=2^N$. 
\begin{lemma} \label{lem:cost}
	\begin{equation}
	C(N,d) \stackrel{.}{=} \frac{N^{d-1}\ 2^{N-d+1}}{(d-1)!},
	\label{eq:effort_leading}
\end{equation}
where $\stackrel{.}{=}$ means the leading order term in $N$.
\end{lemma}
\begin{proof}
	By induction on $d$. The claim holds clearly for $d=1$. Assume, it holds for $d-1$. By considering the recurrence formula \eqref{eq:rekurs_ges}, we see that $C(N,d)=p(N)\ 2^{N}$, where $p$ is a polynomial of order less or equal to $d$. Consequentely,
	\begin{eqnarray*}
		C(N,d) & \stackrel{.}{=} & \frac{N^{d-2}\ 2^{N-d+2}}{(d-2)!} + \sum_{k=1}^{N}\frac{(N-k)^{d-2}\ 2^{N-k-d+2}}{(d-2)!}2^{k-1} \\
					 & = & \frac{N^{d-2}\ 2^{N-d+2}}{(d-2)!} + \frac{2^{N-d+1}}{(d-2)!} \sum_{k=1}^{N}(N-k)^{d-2} \\
					 & \stackrel{.}{=} & \frac{N^{d-2}\ 2^{N-d+2}}{(d-2)!} + \frac{2^{N-d+1}}{(d-2)!}\frac{N^{d-1}}{d-1} \\
					 & \stackrel{.}{=} & \frac{N^{d-1}\ 2^{N-d+1}}{(d-1)!}
	\end{eqnarray*}
\end{proof}

According to (\ref{eq:accuracy2}), the approximation error $\left\|f-f_{\vek{I}}\right\|$ is bounded by $\sum_{\ell\notin\vek{I}}\left\|f_{\ell}\right\|$, i.e.
\[
\left\|f-f_{\vek{I}}\right\| \leq \sum_{|\ell|> N}\left\|f_{\ell}\right\|,
\]
if we use the optimal approximation space $W_{\bI(N)}$. By (\ref{eq:L1beitrag}) this means
\[
\left\|f-f_{\vek{I}}\right\| \leq \sum_{|\ell|> N}\left[2^{-\sum_{\ell_i\neq0}(\ell_i+1)}\prod_{\ell_{i}\neq0} \|\partial_i f\|_\infty\right]
\]
Again, the constants $\prod_{\ell_{i}\neq0} \|\partial_i f\|_\infty$ only depend on the function to be approximated.  Thus, without a priori knowledge about $f$ we need to assume that they can be bounded by some common constant and accordingly define the discretization error of the $N^{th}$ level sparse basis as
\begin{equation}
	E(N,d) = \sum_{|\ell|>N} 2^{-\sum_{\ell_i\neq0}(\ell_i+1)}.
	\label{eq:sparse_error}
\end{equation}
Let $E(-n,d)$ for $n\in\nat, n>0$ represent the error of the empty basis and $\ell = (\tilde{\ell},\ell_d)$ with $\tilde{\ell}\in\nat_{0}^{d-1}$. Then 
\begin{eqnarray*}
	E(N,d) & = & \sum_{|\ell|>N} 2^{-\sum_{\ell_i\neq0}(\ell_i+1)} \\
				 & = & \sum_{\ell_d=0}^{\infty}2^{-(\ell_d+1)(\ell_d\neq0)} \sum_{|\tilde{\ell}|>N-\ell_d}2^{-\sum_{\tilde{\ell}_i\neq0}(\tilde{\ell}_i+1)} \\
				 & = & \sum_{\ell_d=0}^{\infty}2^{-(\ell_d+1)(\ell_d\neq0)} E(N-\ell_d,d-1),
\end{eqnarray*}
where the expression $(\ell_i\neq0)$ has the value $1$, if it is true, otherwise $0$. This leads, by splitting the sum, to the recurrence formula
\begin{equation}
	E(N,d) = E(N,d-1) + \sum_{k=1}^{N}E(N-k,d-1)2^{-k-1} + \underbrace{\sum_{k=N+1}^{\infty}2^{-k-1}E(-1,d-1)}_{=2^{-N-1}E(-1,d-1)}.
	\label{eq:rekurs_err}
\end{equation}
We easily compute that $E(N,1) = 2^{-N-1}$ for $N\geq0$ and $E(-1,d)=(3/2)^d$. 
\begin{lemma}\label{lem:error}
	\begin{equation}
	E(N,d) \stackrel{.}{=} \frac{N^{d-1}\ 2^{-N-d}}{(d-1)!},
	\label{eq:err_leading}
\end{equation}
where, again, $\stackrel{.}{=}$ means the leading order term in $N$.
\end{lemma}
\begin{proof}
	By induction on $d$. The claim holds for $d=1$, assume it holds for $d-1$. Then
\begin{eqnarray*}
	E(N,d) & \stackrel{.}{=} & \frac{N^{d-2}2^{-N-d+1}}{(d-2)!} + \sum_{k=1}^{N}\frac{(N-k)^{d-2}2^{-N+k-d+1}}{(d-2)!}2^{-k-1} + \left(\frac{3}{2}\right)^{d-1}2^{-N-1} \\
				 & \stackrel{.}{=} & \frac{N^{d-2}2^{-N-d+1}}{(d-2)!} + \frac{2^{-N-d}}{(d-2)!}\sum_{k=1}^{N}(N-k)^{d-2} \\
				 & \stackrel{.}{=} & \frac{2^{-N-d}}{(d-2)!}\frac{N^{d-1}}{d-1}
\end{eqnarray*}
\end{proof}

\paragraph{Comparison with Ulam's method.}

We now compare the expressions for the asymptotic behaviour of cost and discretization error in dependence of the discretization level $N$ and the problem dimension $d$ in Lemmata~\ref{lem:cost} and \ref{lem:error} to the corresponding expressions for the standard Ulam basis, i.e.\ the span of the characteristic functions on a uniform partition of the unit cube into cubes of edge length $2^{-M}$ in each coordinate direction -- this is $\bigoplus_{\norm{\ell}{\infty}\leq M}W_\ell$. This space consists of $(2^M)^d$ basis functions, the discretization error is $\mathcal{O}\left(2^{-M  }\right)$.

We thus have -- up to constants -- the following asymptotic expressions for cost and error of the sparse and the standard basis:
\begin{center}
	\begin{tabular}{c|c|c}
		 & cost & error \\ \hline
		sparse basis & $ (N/2)^{d-1}\, 2^{N}$ & $(N/2)^{d-1}\, 2^{-N}$ \\ \hline
		standard basis & $2^{dM}$ & $2^{-M}$
	\end{tabular}
\end{center}

To highlight the main difference, consider the following simple computation: The expressions for the errors are equal if
\[ 
M = N+d-(d-1)\log_{2}N. 
\]
Using this value for $M$ in the cost expression we get 
$N^{d-1}\ 2^{N-d} < 2^{dN+d^2-d(d-1)\log_{2}N}$, i.e.
\begin{equation}
	\frac{N}{d+1} > \log_{2}N-1
	\label{eq:sparseVSstandard}
\end{equation}
as a sufficient condition for the sparse basis to be more efficient than the standard basis. Since we neglected constants and lower order terms in this estimate, the only conclusion we can draw from this is that from a certain accuracy requirement on, the sparse basis is more efficient than the standard one.

\subsection{Computing the matrix entries}

When we use Monte-Carlo quadrature in order to approximate the entries of the transition matrix in both methods, the overall computation breaks down into the following three steps:
\begin{enumerate}
	\item mapping the sample points,
	\item constructing the transition matrix,
	\item solving the eigenproblem.
\end{enumerate}
While steps 1.\ and 3.\ are identical for both methods, step 2.\ differs significantly.  This is due to the fact that in contrast to Ulam's method, the basis functions of the sparse hierarchical tensor basis have global and non-disjoint supports. 

\subsubsection{Number of sample points}\label{sssec:sample}

Applying Monte-Carlo approximation to  \eqref{eq:FPO_sparse2}, we obtain
\begin{eqnarray}
	\tilde{p}_{ij} &=& |\varphi_i||\varphi_j|\left( \frac{m\left(X_j^{+}\right)}{K_j}\sum_{k=1}^{K_j}\chi_{i}^+\left(S(x_k^+)\right) - \chi_{i}^-\left(S(x_k^+)\right)\right.\\
	&-& \left.\frac{m\left(X_j^{-}\right)}{K_j}\sum_{k=1}^{K_j}\chi_{i}^+(S(x_k^-)) - \chi_{i}^-(S(x_k^-))\right),
	\label{eq:computePij2}
\end{eqnarray}
where the sample points $x_k^\pm$ are chosen i.i.d.\ from a uniform distribution on $X_j^+$ and $X_j^-$, respectively.  In fact, since the union of the supports of the basis functions in one subspace $W_\ell$ covers all of $X$, we can reuse the same set of $\varrho$ sample points and their images for each of the subspaces $W_\ell$ (i.e.\ ${N+d}\choose{d}$ times).  Note that the number $K_j$ of test points chosen in $X_j^\pm$ now varies with $j$ since the supports of the various basis functions are of different size: on average, $K_j=\varrho m(X_j^\pm)$. Accordingly, for the absolute error of $\tilde p_{ij}$ we get
\begin{equation}
	|\tilde{p}_{ij}-p_{ij}| \sim \frac{m(X_{j})}{\sqrt{m(X_{i})m(X_{j})}\sqrt{\varrho\ m(X_{j})}} = \frac{1}{\sqrt{\varrho m(X_i)}},
	\label{eq:entr_err_spu}
\end{equation}
where we used that $m(X_{i}^{\pm}) \sim m(X_{i})$. In the worst case we thus get
\[	|\tilde{p}_{ij}-p_{ij}| \sim \frac{2^{N/2}}{\sqrt{\varrho}}, \]
which implies 
\begin{equation}
	\varrho \apprge \frac{2^N}{TOL^2}.
	\label{eq:point_nr_spu}
\end{equation}
for the total number of test points required in order to achieve an accuracy of $TOL$ in the entries of the transition matrix.

\paragraph{Comparison with Ulam's method.}

Aiming at a final accuracy of $\epsi>0$ of the eigenvector, we have to choose $M$ and $N$ accordingly.  Assuming that the corresponding eigenproblems are well conditioned, $TOL = \epsi$ is a reasonable choice for the required accuracy of the entries. This implies a number of
\[ \varrho \apprge \epsi^{-(d+2)} \]
sample points for the standard realisation of Ulam's method (cf. \ref{ssec:comp_trm}), and yields, since $2^N \apprle \epsi^{-1}\left(\log(\epsi^{-1})\right)^{d-1}$, 
\[ \varrho \apprge \epsi^{-3}\left(\log(\epsi^{-1})\right)^{d-1} \]
sample points for the sparse Ulam method.  Note that for $d\geq 2$, the sparse Ulam method requires less sample points than Ulam's method in order to achieve a comparable accuracy in the eigenvector approximation.

\subsubsection{Number of index computations}\label{sssec:index}

While in Ulam's method each sample point is used in the computation of one entry of the transition matrix only, this is not the case in the Sparse Ulam method.  In fact, each sample point (and its image) is used in the computation of $|\bI(N)|^2$ matrix entries, namely one entry for each pair $(W_\ell,W_{m})$ of subspaces.  

Correspondingly, for each sample point $x$ (and its image) and for each $\ell\in \bI(N)$, we have to compute the index $\vj$ of the basis function $\varphi_{\ell,\vj}\in W_\ell$ whose support contains $x$.  Since (cf.\ the previous section) the required number of sample points is $\mathcal{O}\left(\frac{2^N}{TOL^2}\right)$ and $|\bI(N)|={N+d \choose d}$, this leads to 
\[ 
\varrho\cdot |\bI(N)| = \frac{2^N}{TOL^2}{N+d \choose d} \lesssim \frac{1}{TOL^2}\frac{N^d}{d!}2^N \stackrel{.}{=} \frac{2^{d-1} N}{d} \frac{1}{TOL^2} \dim V_N
\]
of these computations (for reasonable $d$).  In contrast, in Ulam's method, the corresponding number is
\[ 
\varrho \cdot 1 = \frac{2^{dM}}{TOL^2} = \frac{1}{TOL^2} \dim V_M.
\]
Note that for the Sparse Ulam method the number of index computations is not staying proportional to the dimension of the approximation space.  However, it is still scaling much more mildly with $d$ than for Ulam's method.

\subsubsection{Occupancy of the transition matrix}

The matrix which represents the discretized transfer operator in Ulam's method is \emph{sparse}: the supports of the basis functions are disjoint, and thus $p_{ij} \neq 0$ only if $S(X_j)\cap X_i \neq \emptyset$.  Hence, for a sufficiently fine partition, the number of partition elements $X_i$ which are intersected by the image $S(X_j)$ is determined by the local expansion of $S$.  This is a fixed number related to a Lipschitz estimate on $S$ and so the matrix of the discretized transfer operator with respect to the standard Ulam basis is sparse for sufficiently large $n$.  Unfortunately this property is not shared by the matrix with respect to the sparse basis as the following considerations show.

\begin{figure}[ht]
	\centering
		\includegraphics[width=0.6\textwidth]{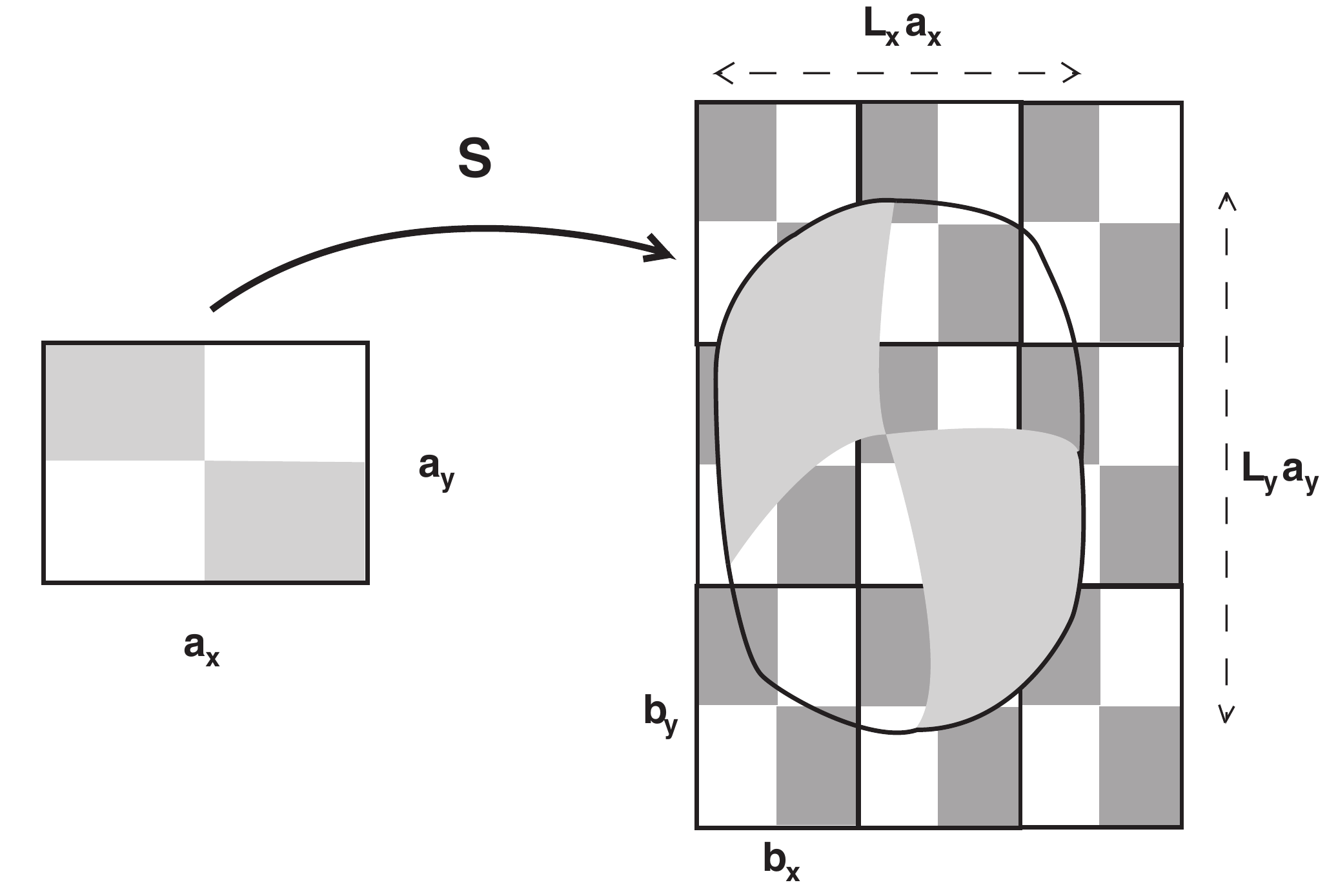}
	\caption{Model for the matrix occupancy in 2D. Shaded and colorless (white) show the function values ($\pm|\varphi|$), thicker black lines the support boundaries.}
	\label{fig:stretch}
\end{figure}

The main reason for this is that the supports of the basis functions in the sparse basis are not localised, cf.\ the thin and long supports of the basis of $W_{\ell}$ for $\ell = (N,0,\ldots,0)$. This means that the occupancy of the transition matrix strongly depends on the \textit{global behaviour} of the dynamical system $S$. 
Let $B_\ell := \left\{ \varphi_{\ell,\vj}\mid \vj_{i}\in \{0,\ldots,2^{\ell_i}-1\}  \right\}$  denote the basis of $W_\ell$ and let
\[ 
\nnz(\bk,\ell) = \left|\left\{(i,j)\mid S(\text{supp}(\varphi_i))\cap\text{supp}(\varphi_j)\neq\emptyset, \varphi_i\in B_{\vek{k}}, \varphi_j\in B_{\ell} \right\}\right| 
\] 
be the number of nonzero matrix entries which arise from the interaction of the basis functions from the subspaces $W_{\vek{k}}$ and $W_\ell$ if $W_{\vek{k}}$ is mapped. 
We define the \emph{matrix occupancy} of a basis $B_\bI=\bigcup_{\ell\in\bI} B_\ell$ as
\begin{equation}
	\nnz(B_\bI) = \sum_{\vek{k},\ell\in \bI} \nnz(\vek{k},\ell).
	\label{eq:occu1}
\end{equation}
In order to estimate $\nnz(\vek{k},\ell)$ we employ upper bounds $L_i$, $i=1,\ldots,d$, for the Lipschitz-constants of $S$, cf.\ Figure~\ref{fig:stretch}. We obtain

\begin{proposition}
\begin{equation}
	\nnz(\vek{k},\ell) \leq \left|B_{\vek{k}}\right|\prod_{i=1}^{d} \left\lceil \frac{L_i\cdot 2^{-\vek{k}_i+1-(\vek{k}_i=0)}} {2^{-\ell_i+1-(\ell_i=0)}}\right\rceil.
	\label{eq:occu2}
\end{equation}
\end{proposition}
\begin{proof}
	Since we have used upper bounds for the Lipschitz constants, \textbf{one} mapped box has at most the extension $L_i\cdot 2^{-\vek{k}_i+1-(\vek{k}_i=0)}$ in the $i^{th}$ dimension. Consequently, its support intersects with at most
	\[ \left\lceil \frac{L_i\cdot 2^{-\vek{k}_i+1-(\vek{k}_i=0)}} {2^{-\ell_i+1-(\ell_i=0)}}\right\rceil \]
	supports of basis functions from $W_{\ell}$.
\end{proof}
\begin{bem}
	Numerical experiments suggest that the above bound approximates the matrix occupancy pretty well. However, it could be improved: \eqref{eq:FPO_sparse1} shows that a matrix entry could still be zero even if supp$(\varphi_i)$ and supp$(P\varphi_j)$ intersect. This is e.g.\ the case if supp$(P\varphi_j)$ is included in a subset of supp$(\varphi_i)$, where $\varphi_i$ is constant (i.e.\ does not change sign). The property $\normLone{Pf}=\normLone{f}$ for $f\geq0$ and positivity (see \cite{LaMa94a}) of $P$ imply $p_{ij}=0$, since $\normLone{\varphi_j^{+}} = \normLone{\varphi_j^{-}}$.
\end{bem}

\paragraph{An asymptotic estimate.}

Let us examine $\nnz(\bk,\ell)$ for $\vek{k} = (0,\ldots,0,N)$ and $\ell = (N,0,\ldots,0)$. By taking all Lipschitz-constants $L_i=1$ we get 
\[ \nnz(\bk,\ell) \apprge  2^{2N}, \]
since $|B_{\vek{k}}|=2^{N-1}$ and the image of each basis function from $B_\bk$ intersects with each basis function from $B_{\ell}$.  Since $|B_N|\approx N^{d-1}2^N$, we get
\begin{equation}
	2^{2N} \apprle \nnz(B_N) \apprle N^{2d-2}2^{2N}.
	\label{eq:occu_bound1}
\end{equation}
The exponential term dominates the polynomial one for large $N$, so asymptotically we will not get a sparse matrix. 

Does this affect the calculations regarding efficiency made above? As already mentioned, the error of Ulam's method is $\epsi=\mathcal{O}(2^{-M})$ while its cost is $2^{dM} = \mathcal{O}(\epsi^{-d})$. Assuming that the Sparse Ulam method has the same error $\epsi = \mathcal{O}(N^{d-1}2^{-N})$, its worst-case cost is 
\[ \mathcal{O}(N^{2d-2}2^{2N}) \apprle \epsi^{-2}N^{4d-4} \apprle \epsi^{-2}\log(\epsi^{-1})^{4d-4}, \]
where we used $N^{d-1}2^{-N} \apprle 2^{-N/2}$, which leads to $\log(\epsi) \apprle -N$. Clearly, this means -- similarily to subsection \ref{subsec:cost_accu} -- partially overcoming the curse of dimensionality.
Even in the most optimistic case, ie. the costs are of $\mathcal{O}(2^{2N})$, we have at least $\mathcal{O}(\epsi^{-2}N^{2d-2})$ costs, so the sparse-Ulam-method is efficienter than Ulam's, only if $d\geq3$.

\section{Numerical examples}\label{sec:examples}

\subsection{A 3d expanding map}
We compare both methods by approximating the invariant density of a simple three dimensonal map. Let $S_i:[0,1]\rightarrow[0,1]$ be given by
\begin{eqnarray*}
 S_1(x) &=& 1-2|x-1/2|, \\
 S_2(x) &=& \left\{ \begin{array}{ll}
	2x/(1-x), & x<1/3\\
	(1-x)/(2x), & \text{else},
	\end{array} \right.,\\
 S_3(x) &=& \left\{ \begin{array}{ll}
	2x/(1-x^2), & x<\sqrt{2}-1\\
	(1-x^2)/(2x), & \text{else},
	\end{array} \right.
\end{eqnarray*}
and $S:[0,1]^3\rightarrow[0,1]^3$ be the tensor product map
\[ S(x) = \left(S_1(x_1),S_2(x_2),S_3(x_3)\right)^{\top}, \]
where $x = (x_1,x_2,x_3)^{\top}$.  This map is expanding and its unique invariant density is given by
\[ 
h(x) = \frac{8}{\pi(1+x_3^2)(1+x_2)^2}. 
\]
(cf.\ \cite{DiZh96a}).

We approximate $h$ by Ulam's method on an equipartition of $2^{3M}$ boxes for $M=4,5,6$ as well as by the Sparse Ulam method on levels $N=4,5,6$.
\begin{figure}[htb]
	\centering
		\includegraphics[width=0.49\textwidth]{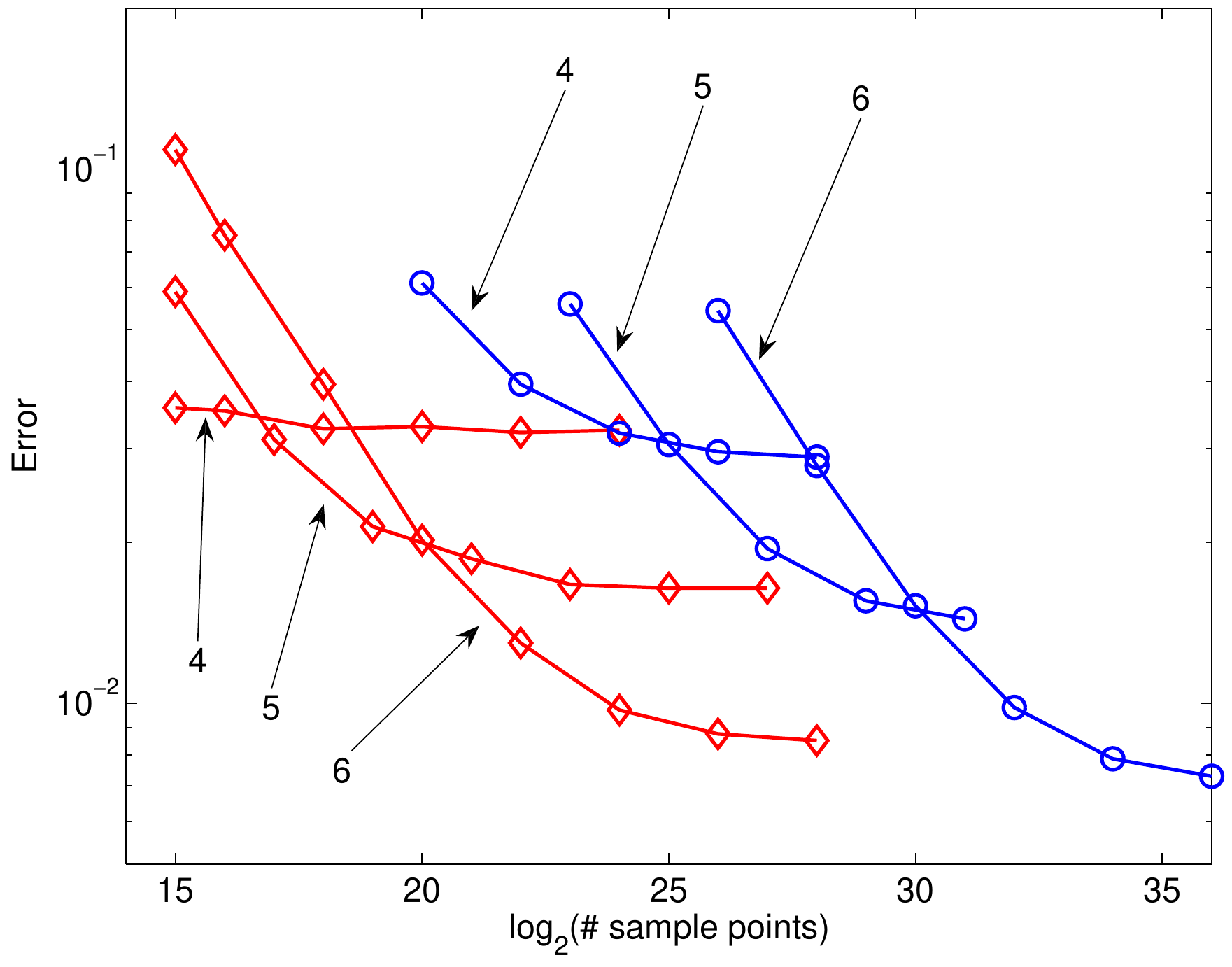}
		\includegraphics[width=0.49\textwidth]{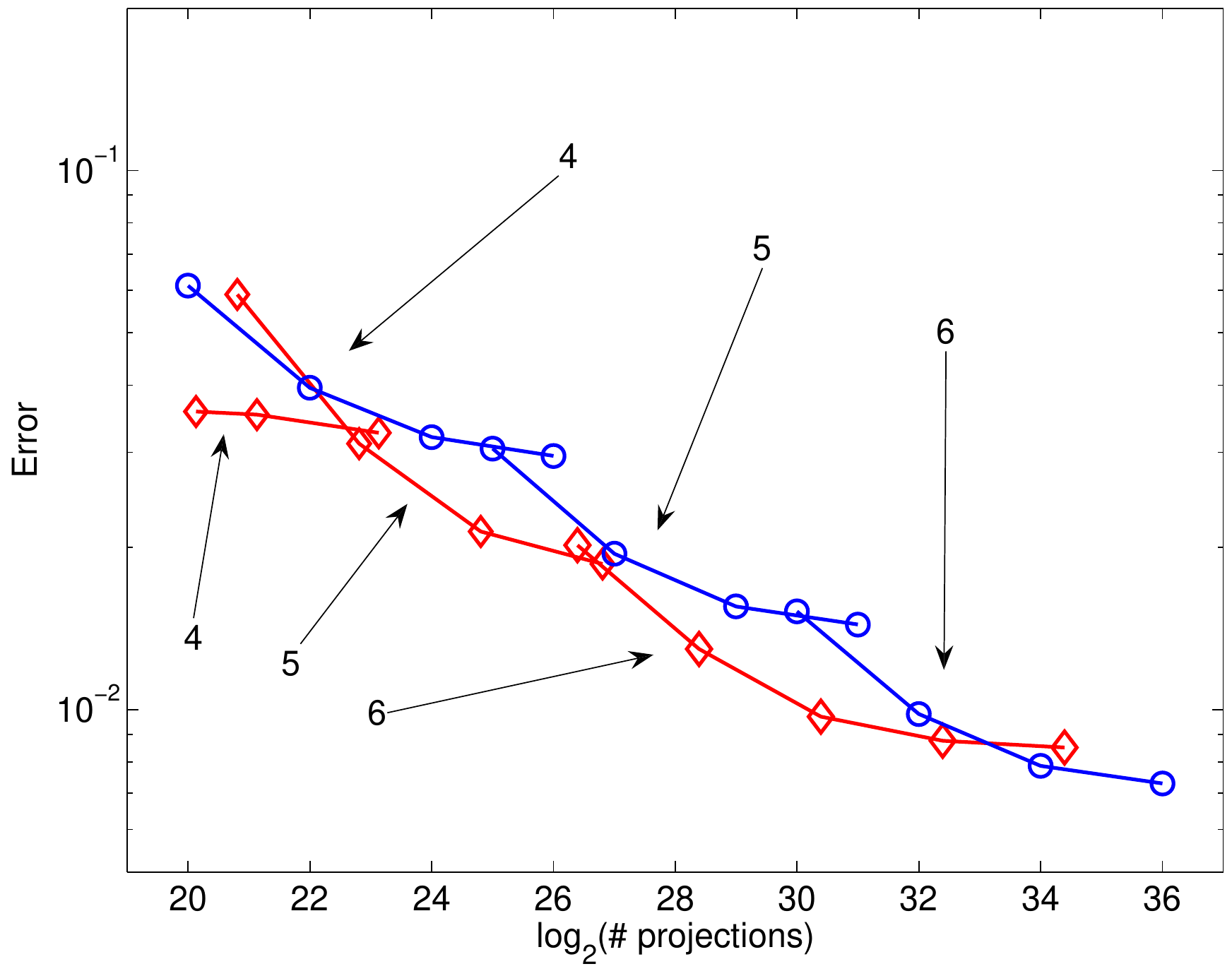}
	\caption{Left: $L^1$-error of the approximate invariant density in dependence on the number of sample points for levels $N,M=4,5,6$. Right: Corresponding number of index computations.}
	\label{fig:Points}
\end{figure}
Figure~\ref{fig:Points} shows the $L^1$-error for both methods in dependence of the number of sample points (left) as well as the number of index computations along these curves (right).  While the Sparse Ulam method requires almost three orders of magnitude fewer sample points than Ulam's method, the number of index computations is roughly comparable.  This is in good agreement with our theoretical considerations in sections~\ref{sssec:sample} and \ref{sssec:index}. 

\begin{figure}[htb]
	\centering
		\includegraphics[width=0.5\textwidth]{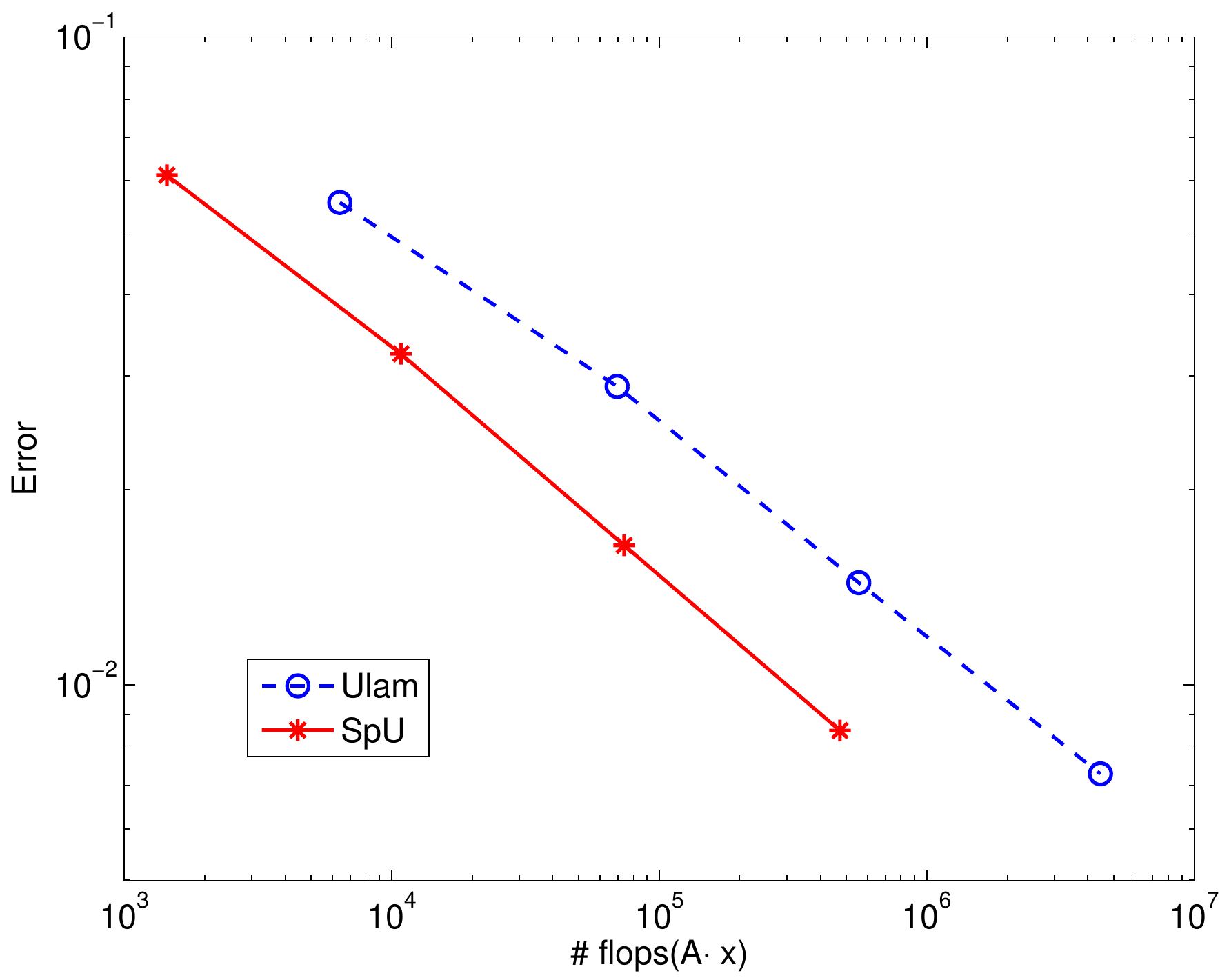}
	\caption{$L^1$-error of the approximate invariant densities in dependence on the number of nonzeros in the transition matrices.}
	\label{fig:Flops}
\end{figure}

In Figure~\ref{fig:Flops} we show the dependence of the $L^1$-error on the number of nonzeros in the transition matrices for levels $M,N=3,\ldots,6$.  Again, the Sparse Ulam method is ahead of Ulam's method by almost an order of magnitude. 

\subsection{A 4d conservative map}

In a second numerical experiment, we approximate a few dominant eigenfunctions of the transfer operator for an area preserving map. Since the information on almost invariant sets does not change \cite{Fr05a} (but the eigenproblem becomes easier to solve) we here consider the symmetrized transition matrix $\frac{1}{2}(P+P^\top)$, cf.\ also \cite{JuMaMe04a}.

\begin{figure}[htb]
	\centering
		\includegraphics[width=0.35\textwidth]{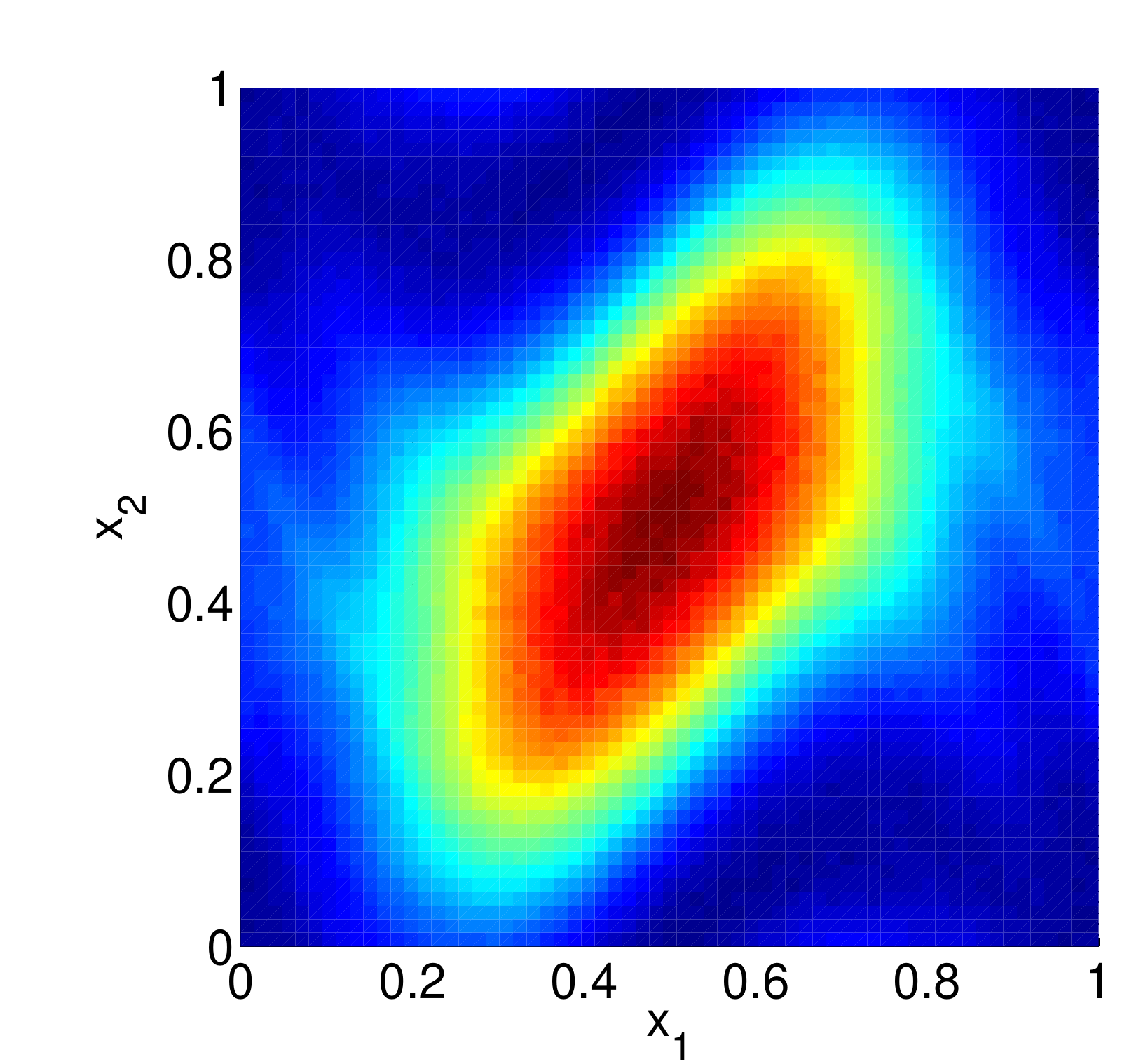}
		\hspace*{2cm}
		\includegraphics[width=0.33\textwidth]{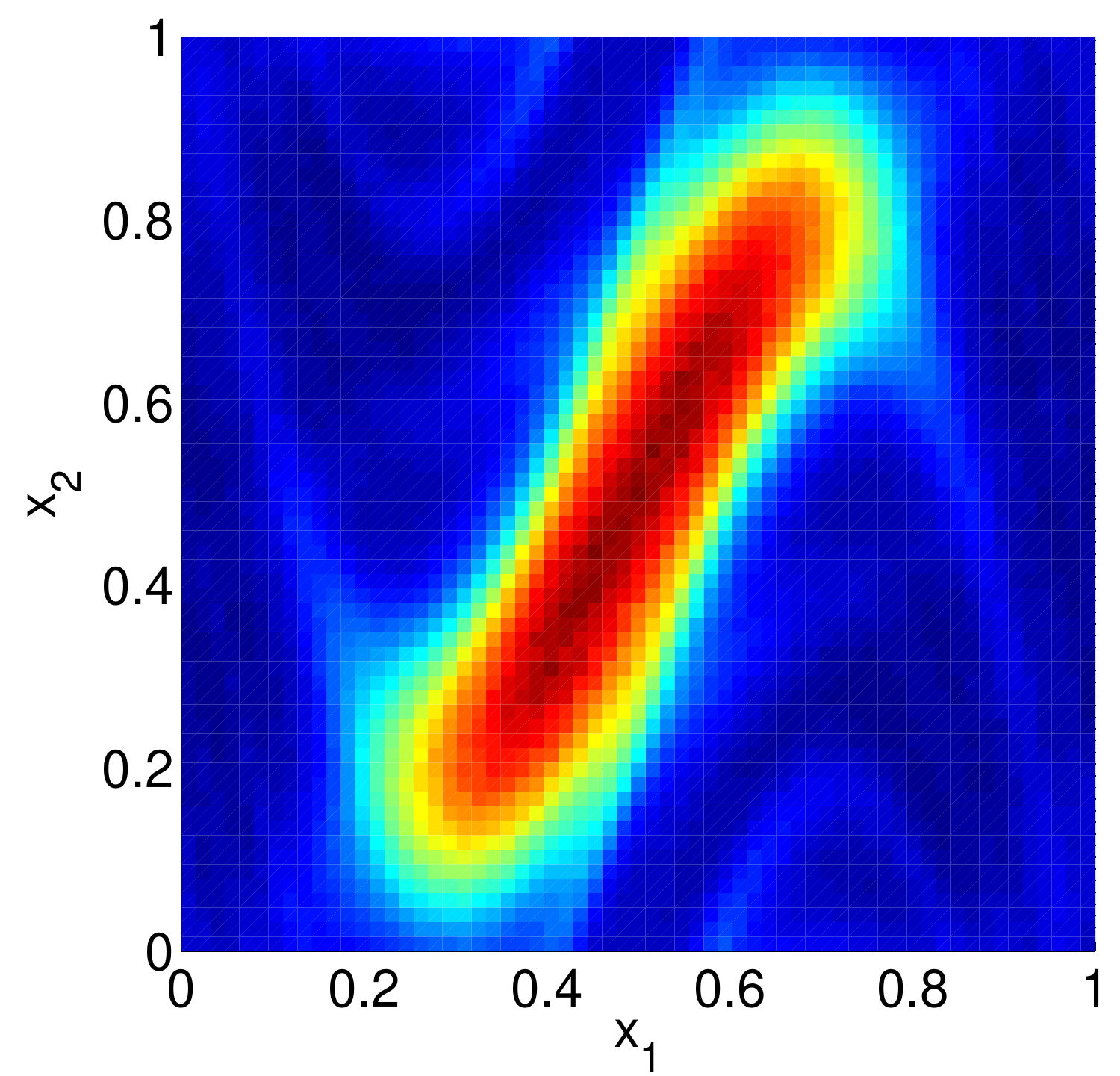}
	\caption{Eigenfunction of the symmetrized transition matrix at the second largest eigenvalue for the standard map. Left: $\rho=0.3$, $\lambda_2=0.97$, right: $\rho=0.6$, $\lambda_2 = 0.93$.}
	\label{fig:std_2d}
\end{figure}

Consider the so called \emph{standard map} $S_{\rho}:[0,1]^2\rightarrow[0,1]^2$,
\[ (x_1,x_2)^{\top} \mapsto \left(x_1+x_2+\rho \sin(2\pi x_1)+0.5, x_2+\rho\sin\left(2\pi x_1\right)\right)^{\top} \mod 1, \]
where $0<\rho<1$ is a parameter.
This map is \emph{area preserving}, i.e.\ the Lebesgue measure is invariant w.r.t. $S_{\rho}$.  Figure~\ref{fig:std_2d} shows approximations of the eigenfunctions at the second largest eigenvalue of $S_\rho$ for $\rho=0.3$ (left) and $\rho=0.6$ (right) computed via Ulam's method on an equipartition of $2^{2\cdot 6}$ boxes (i.e.\ for $M=6$).

\begin{figure}[h]
	\centering
		\includegraphics[width=0.70\textwidth]{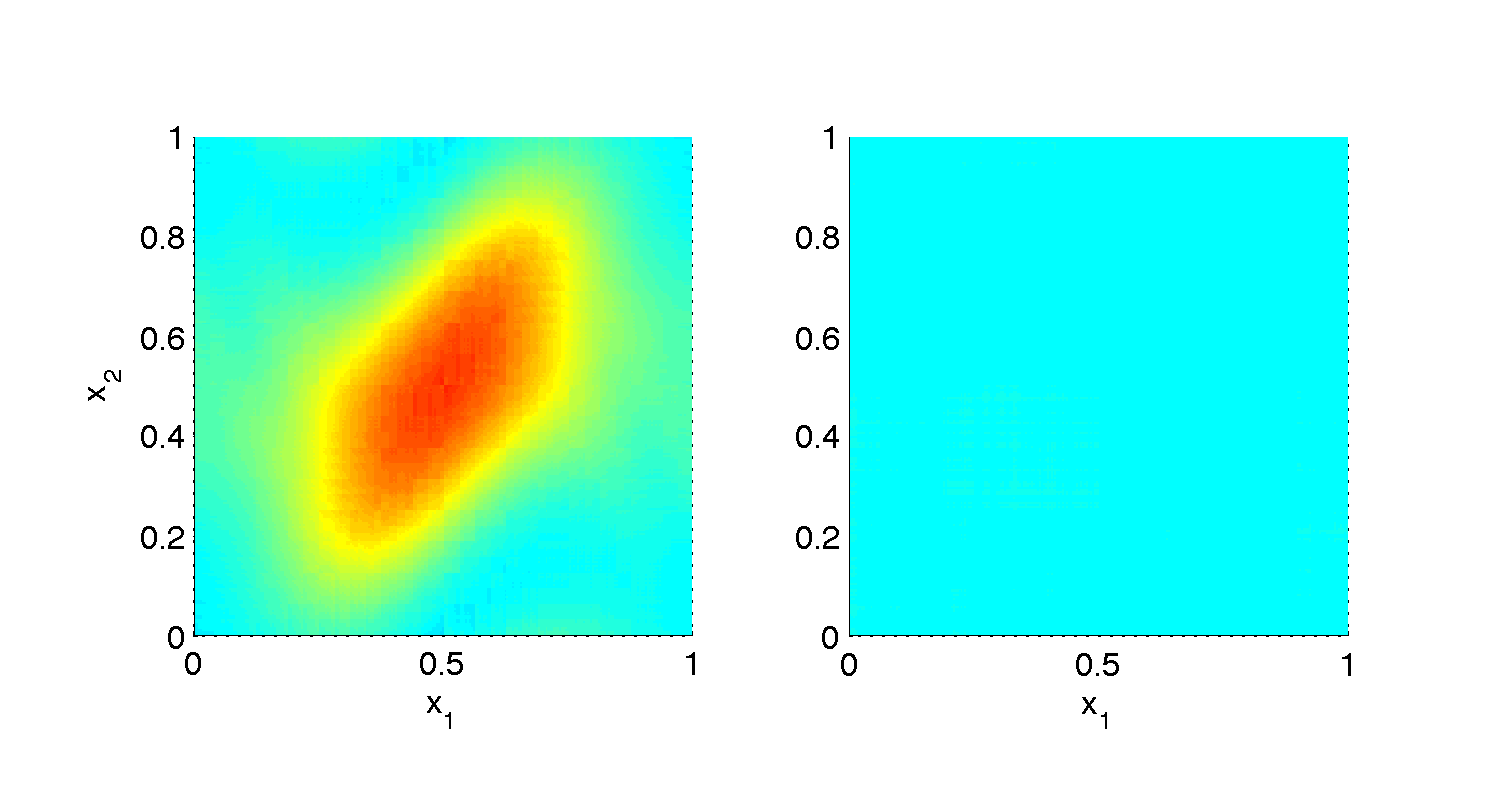}
		\caption{Approximate eigenfunction at $\lambda_2 = 0.97$. Left: $v_2(\cdot,\cdot,x_3,x_4)$ for fixed $x_3,x_4$, right: $v_2(x_1,x_2,\cdot,\cdot)$ for fixed $x_1,x_2$.} 
	\label{fig:std_4d-1}
\end{figure}

\begin{figure}[h]
	\centering
		\includegraphics[width=0.70\textwidth]{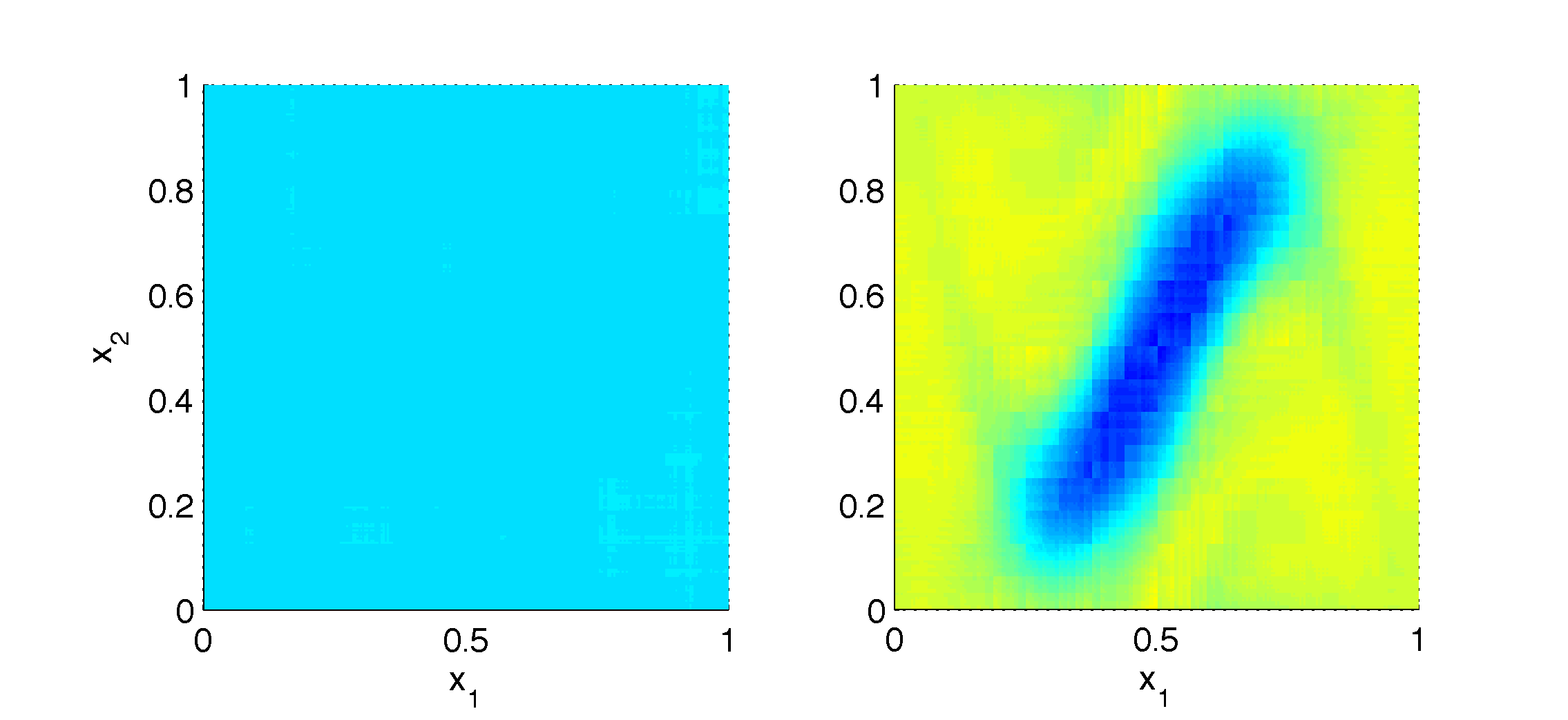}
		\caption{Approximate eigenfunction at $\lambda_2 = 0.93$. } 
	\label{fig:std_4d-2}
\end{figure}

We now define $S:[0,1]^4\rightarrow[0,1]^4$ by
\[
S = S_{\rho_1}\otimes S_{\rho_2},
\]
with $\rho_1 = 0.3$ and $\rho_2 = 0.6$.   Note that the eigenfunctions of $S$ are tensor products of the eigenfunctions of the $S_{\rho_i}$. This is reflected in Figures~\ref{fig:std_4d-1} and \ref{fig:std_4d-2} where we show the eigenfunctions at the two largest eigenvalues, computed by the Sparse Ulam method on level $N=8$, using $2^{24}$ sample points overall. Clearly, each of these two is a tensor product of the (2d-) eigenfunction  at the second largest eigenvalue with the (2d-) invariant (i.e.\ constant) density.

\begin{figure}[h]
	\centering
		\includegraphics[width=0.70\textwidth]{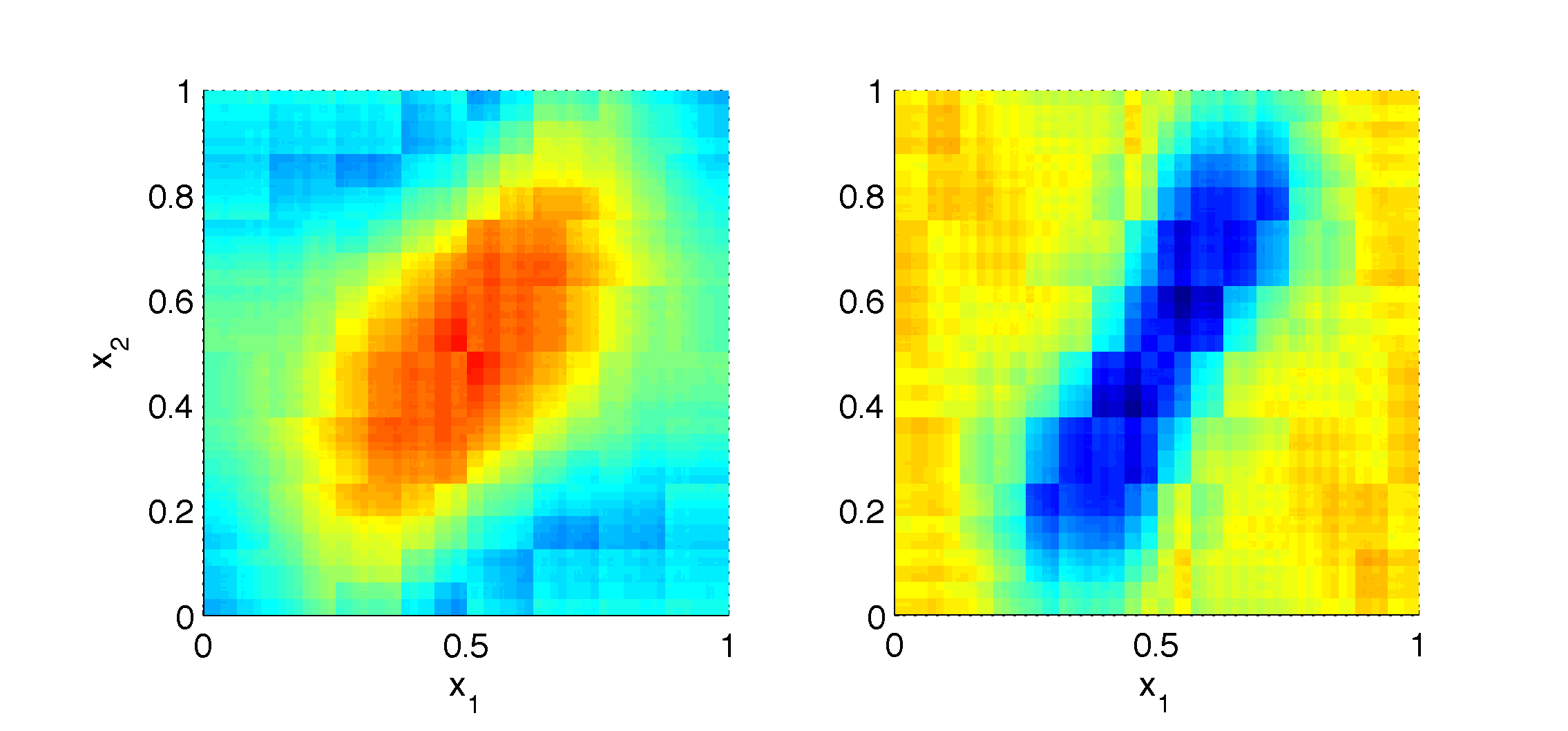}
		\caption{Approximate eigenfunction at $\lambda = 0.80$.} 
	\label{fig:std_4d-3}
\end{figure}

Figure~\ref{fig:std_4d-3} shows an eigenfunction for which both factors of the tensor product are non-constant. The resolution of this eigenfunction seems worse than for those with one constant factor.  In fact, for an approximation of an eigenfunction which is constant with respect to, say, $x_3$ and $x_4$ it suffices to consider subspaces $W_\ell$ with $\ell=(\ell_1,\ell_2,0,0)$. All  other coefficients are zero, the problem reduces to a two-dimensional one and so the eigenfunctions are not perturbed by basis functions varying in the $x_3$ and $x_4$ directions.

\small
\bibliography{bibdesk}
\bibliographystyle{abbrv}

\end{document}